\documentclass[12pt]{amsart}
\usepackage{latexsym, amsmath, amscd, amssymb, amsthm} 
\usepackage[T2A]{fontenc}
 \usepackage[cp1251]{inputenc}
\usepackage[english]{babel}
\usepackage{color}
 \def\mathbi#1{\textbf{\em #1}}
\begin{document}

 \title[]{ The MAPLE package for calculating  Poincar\'e series.}

\author{Leonid Bedratyuk}  
\begin{abstract} 
  We  offer a Maple package {\tt Poincare\_Series}  for calculating the Poincar\'e series for  the algebras of invariants/covariants of binary forms, for the algebras of joint invariants/covariants of several binary forms, for  the kernel of Weitzenb\"ock derivations,for the bivariate Poincar\'e series of algebra of covariants of binary $d$-form and  for  the multivariate Poincar\'e series of the algebras of joint invariants/covariants of several binary forms.
\end{abstract}
\maketitle

\section{ Introduction}

The Poincar\'e  series of a graded algebra $A=\oplus_{i} (A)_i$  is defined as  formal power series  $\mathcal{P}(A,z):=\sum_{i=0}^{\infty} \dim (A)_i z^i.$ If an algebra is finitely generated then its   Poincar\'e series is    the power series expansions  of certain  rational functions.

 The present package implements results of the following papers:

\begin{itemize}
	\item Leonid Bedratyuk, The Poincar\'e  series of the covariants of binary forms, Int. Journal of Algebra,2010

\item Leonid Bedratyuk,The Poincar\'e series of the joint invariants and covariants of the two binary forms,  Linear and Multilinear algebra, 2010

\item Leonid Bedratyuk, Linear locally nilpotent derivations and the classical invariant theory, I: The Poincare series, Serdica Math. J.,2010

\item Leonid Bedratyuk, Bivariate Poincar\'e series for the algebra of covariants of a binary form, preprint arXiv:1006.1974 

\item Leonid Bedratyuk, Multivariate Poincar\'e series for the algebras of joint  invariants and covariants of several binary forms, in preparation.
\end{itemize}

\section{ Instalation.}  The package can be downloaded  from the web  http://sites.google.com/site/bedratyuklp/.
 
\begin{enumerate}
\item download the file {\tt Poincare\_Series.mpl} and save it into your  Maple directory;
\item download the Xin's file (see a link at the web page) {\tt Ell2.mpl} and save it into your  Maple directory;
 \item  run Maple;
 
 \item  \textcolor{red}{{\tt > read "Poincare\_Series.mpl": read "Ell2.mpl":}}

 \item  If necessary use \textcolor{red}{{\tt > Help();}}
\end{enumerate}

\section{ Formulas for the Poincar\'e  series. }

 Below are  the list of main formulas.

\subsection{ Invariants and covariants of   binary form}  Let $\mathcal{I}_d,$ $\mathcal{C}_d $ be algebras of invariants and covariants of binary $d$-form  graded under degree.   We  have 
\begin{gather}
\displaystyle \mathcal{P}(\mathcal{I}_d,z)=\sum_{0\leq k <d/2} \varphi_{d-2\,k} \left( \frac{(-1)^k z^{k(k+1)} (1-z^2)}{(z^2,z^2)_k\,(z^2,z^2)_{d-k}} \right), \text{ (Springer's formula)},\\
\displaystyle \mathcal{P}(\mathcal{C}_d,z)=\sum_{0\leq k <d/2} \varphi_{d-2\,k} \left( \frac{(-1)^k z^{k(k+1)} (1+z)}{(z^2,z^2)_k\,(z^2,z^2)_{d-k}} \right),
\end{gather}
here $(a,q)_n=(1-a) (1-a\,q)\cdots (1-a\,q^{n-1})$ denotes the $q$-shifted factorial and the function $\varphi_{n}:\mathbb{C}[[z]] \to \mathbb{C}[[z]]$  defined by 
$$
\varphi_{n}\left(\sum_{i=0}^{\infty} a_{i}z^i \right)=\sum_{i=0}^{\infty} a_{i n} z^i. 
$$
\subsection{ Joint invariants and covariants of   binary form}
Let $\mathcal{I}_{\mathbi{d}},$ $\mathcal{C}_{\mathbi{d}}, $ $\mathbi{d}=(d_1,d_2,\ldots,d_n)$ be algebras of joint invariants and joint covariants of  $n$ binary forms of degrees $d_1,$  $d_2,$ $\ldots,d_n.$  Then 

\begin{gather}
\displaystyle \mathcal{P}(\mathcal{I}_{{\mathbi{d}}},z)=\sum_{i=0}^{d^*} \sum_{k=1}^{\beta_i} \frac{1}{(k-1)!} \frac{d^{k-1}\left(  z^{k-1} \varphi_{d^*-k}((1-z^2)\,A_{i,k}(z))\right)}{dz^{k-1}},\\
\displaystyle \mathcal{P}(\mathcal{C}_{{\mathbi{d}}},z)=\sum_{i=0}^{d^*} \sum_{k=1}^{\beta_i} \frac{1}{(k-1)!} \frac{d^{k-1}\left(  z^{k-1} \varphi_{d^*-k}((1+z)\,A_{i,k}(z))\right)}{dz^{k-1}},
\end{gather}

$$
\begin{array}{l}
\displaystyle  A_{i,k}(z)=\frac{(-1)^{\beta_i-k}}{(\beta_i-k)!\,(z^{i})^{{\beta_i-k}}}\lim_{t \to z^{-i}} \frac{\partial^{{\beta_i-k}}}{\partial t^{{\beta_i-k}}}\left(f_{d}(tz^{d^*},z)(1-t z^{i})^{\beta_i} \right).  
\end{array}
$$
 The integer numbers $\beta_i,i=0,\ldots,2d^*,$  $d^*:=\max(d_1,d_2,\ldots, d_n),$ are defined  from the decomposition 
 $$
f_{\mathbi{d}}(tz^{d^*},z)=\left((1-t)^{\beta_0} (1-t z)^{\beta_1} (1-t z^2)^{\beta_2} \ldots (1-t z^{2\,d^*})^{\beta_{2\,d^*}}\right)^{-1}, 
$$
where 
$$
f_{\mathbi{d}}(t,z)=\left(\displaystyle \prod_{k=1}^s (tz^{- \,d_k},z^2)_{d_k+1}\right)^{-1}.
$$

\subsection{ Joint invariants and covariants of   linear and quadratic binary forms}
 Let   $d_1=d_2=\ldots=d_n=1,$  i.e. $\mathbi{d}=(1,1,\ldots,1).$ Then 
\begin{gather} 
 \mathcal{P}(\mathcal{I}_{{\mathbi{d}}},z)=\sum_{k=1}^n \frac{(-1)^{n-k}}{(k-1)!} \frac{(n)_{n-k}}{(n-k)! } \frac{d^{k-1}}{dz^{k-1}}\left( \left( \frac{z}{1-z^2} \right)^{2n-k-1} \right)
=\frac{{\rm N}_{n-2}(z^2)}{(1-z^2)^{2n-3}},
\end{gather}
where  ${\rm N}_{n}(z)$ is the $n$-th Narayama polynomial 
$$
{\rm N}_{n}(z)=\sum_{k=1}^n \frac{1}{k} { n-1 \choose k-1}{ n \choose k-1} z^{k-1}. 
$$
and 
$(n)_m:=n(n+1)\cdots (n+m-1),$ $(n)_0:=1$ denotes the shifted factorial.
\begin{gather}
\displaystyle 
\displaystyle \mathcal{P}(\mathcal{C}_{{\mathbi{d}}},z)=\sum_{k=1}^n \frac{(-1)^{n-k}}{(k-1)!} \frac{(n)_{n-k}}{(n-k)! } \frac{d^{k-1}}{dz^{k-1}} \left( \frac{(1+z)z^{2n-k-1}}{(1-z^2)^{2n-k}} \right),
\end{gather}

Let  $d_1=d_2=\ldots=d_n=2,$ ${\mathbi{d}}=(2,2,\ldots,2)$, then 
\begin{gather}
\displaystyle \mathcal{P}(\mathcal{I}_{{\mathbi{d}}},z){=}\sum_{k=1}^n \frac{(-1)^{n{-}k}}{(n{-}k)(k{-}1)!}\frac{d^{k{-}1}}{dz^{k{-}1}} \left(\sum_{i=0}^{n-k} { n{-}k  \choose i}  \frac{(n)_i (n)_{n{-}k{-}i}(1{-}z) z^{2n-k-i-1}}{(1{-} z)^{n+i} (1{-}z^2)^{2n-k-i}}\right),\\
\displaystyle \mathcal{P}(\mathcal{C}_{{\mathbi{d}}},z)=\sum_{k=1}^n \frac{(-1)^{n-k}}{(n-k)!(k-1)!}\frac{d^{k-1}}{dz^{k-1}} \left( \sum_{i=0}^{n-k} { n-k  \choose i}  \frac{(n)_i (n)_{n-k-i} z^{2n-k-i-1}}{(1- z)^{n+i} (1-z^2)^{2n-k-i}}\right)=\\
\displaystyle =\frac{\displaystyle \sum_{i=0}^{n-1} {n-1 \choose i}^2 (z^2)^i}{\displaystyle(1-z)^n(1-z^2)^{2n-1}}.
\end{gather}

\subsection{ Kernel of Weitzenb\"ok derivation}
 Denote by  $\mathcal{D}_{\mathbi{d}}$    the   Weitzenb\"ok derivation (linear locally nilpotent derivation) with   its matrix consisting  of  $n$  Jordan blocks of size  $d_1+1,$ $d_2+1,$ $\ldots, d_s+1,$ respectively. Since $\ker \mathcal{D}_{\mathbi{d}} \cong \mathcal{C}_{{\mathbi{d}}}$ and  the isomorphism preserve  degrees then   have  that  $\mathcal{P}(\ker \mathcal{D}_{\mathbi{d}},z)=\mathcal{P}(\mathcal{C}_{{\mathbi{d}}},z).$

\subsection{ Bivariate Poincare series for covariants of binary form}
The algebra  $\mathcal{C}_{d}$ of covariants is a finitely  generated bigraded algebra: 
$$
\mathcal{C}_{d}=(\mathcal{C}_{d})_{0,0}+(\mathcal{C}_{d})_{1,0}+\cdots+(\mathcal{C}_{d})_{i,j}+ \cdots,
$$
where  each subspace   $(\mathcal{C}_{d})_{i,j}$ of covariants of degree $i$  and order $j$ is   finite-dimensional.  We have 

\begin{gather}
\displaystyle \mathcal{P}(\mathcal{C}_{d},z,t)=\sum_{i=0}^{\infty} (\mathcal{C}_{d})_{i,j} z^i t^j = \sum_{0\leq k <d/2} \psi_{d-2\,k} \left( \frac{(-1)^k t^{k(k+1)} (1-t^2)}{(t^2,t^2)_k\,(t^2,t^2)_{d-k}} \right) \frac{1}{1-z t^{d-2\,k}},
\end{gather}
where $\psi_n:\mathbb{Z}[[t]] \to \mathbb{Z}[[t,z]], n \in \mathbb{Z}_{+}$ be  a  $\mathbb{C}$-linear function   defined by  
$$
\psi_{n}\left(t^m\right):=\left \{ \begin{array}{l} z^i t^j, \text{  if  } m=n\,i-j, j<n,  \\ 0, \text{ otherwise. }  \end{array} \right. 
$$
Note that $\mathcal{P}(\mathcal{C}_{d},z,0)=\mathcal{P}(\mathcal{I}_{d},z)$  and $\mathcal{P}(\mathcal{C}_{d},z,1)=\mathcal{P}(\mathcal{C}_{d},z).$

\subsection{Multivariate Poincar\'e series}
The algebra  $\mathcal{C}_{\mathbi{d}}$ is a finitely  generated multigraded algebra under the multidegree-order: 
$$
\mathcal{C}_{\mathbi{d}}=(\mathcal{C}_{\mathbi{d}})_{\mathbi{m},0}+(\mathcal{C}_{\mathbi{d}})_{\mathbi{m},1}+\cdots+(\mathcal{C}_{\mathbi{d}})_{\mathbi{m},j}+ \cdots,
$$
where  each subspace   $(\mathcal{C}_{\mathbi{d}})_{\mathbi{d},j}$ of covariants of multidegree $\mathbi{m}:=(m_1,m_2,\ldots,m_n)$  and order $j$ is   finite-dimensional. The formal power series 
$$
\mathcal{P}(\mathcal{C}_{\mathbi{d}},z_1,z_2,\ldots,z_n,t)=\sum_{\mathbi{m},j=0}^{\infty }\dim((\mathcal{C}_{\mathbi{d}})_{\mathbi{m},j}) z_1^{m_1} z_2^{m_2}\cdots z_n^{m_n} t^j,
$$ 
is called the multivariariate Poincar\'e series   of the algebra of   join covariants  $\mathcal{C}_{\mathbi{d}}.$

The following formula holds:

$$
\mathcal{P}(\mathcal{C}_{\mathbi{d}},z_1,z_2,\ldots,z_n,t)=\underset{\scriptscriptstyle \geq 0}{\Omega}f_{\mathbi{d}}\left(z_1 (t\lambda)^{d_1},z_2 (t\lambda)^{d_2},\ldots,z_s (t\lambda)^{d_s},\dfrac{1}{t\lambda}\right),
$$
where $$
\begin{array}{l}
f_{{\mathbi{d}}}(z_1,z_2,\ldots,z_n,t)=\displaystyle \frac{1}{\displaystyle \prod_{k=1}^n \prod_{j=0}^{d_k} (1-z_k t^{d_k-2\,j})},\\
\end{array}
$$

For the multivariariate Poincar\'e series   of the algebra of   join invariants  $\mathcal{I}_{\mathbi{d}}$  we   have 
$$
\mathcal{P}(\mathcal{I}_{\mathbi{d}},z_1,z_2,\ldots,z_n)=\underset{\scriptscriptstyle = 0}{\Omega}f_{\mathbi{d}}\left(z_1 (t\lambda)^{d_1},z_2 (t\lambda)^{d_2},\ldots,z_s (t\lambda)^{d_s},\dfrac{1}{t\lambda}\right).
$$
Here  $\underset{\scriptscriptstyle \geq 0}{\Omega}$ and $\underset{\scriptscriptstyle = 0}{\Omega}$  are  the  MacMahon's  Omega operators.

 
\section{Package Commands and Syntax}


{\noindent
{\em Command name}: {\tt INVARIANTS\_SERIES}\\
\noindent{\em Feature}: Computes the Poincare series for the algebras of   joint invariants for the binary forms of degrees $d_1, d_2, \ldots,d_n.$\\
\noindent{\em Calling sequence}:
{\tt
INVARIANTS\_SERIES$\left([d_1, d_2, \ldots,d_n] \right);$
}\\
\noindent{\em Parameters}:\\
{\raggedright
\begin{tabular}{lcl}
$ [d_1, d_2, \ldots,d_n]$ & - & a list of degrees of $n$  binary forms.\\
$ n $ & - & an integer, $n\geq 1.$\\

\end{tabular}\\
\vspace{5mm}
{\noindent
{\em Command name}: {\tt COVARIANTS\_SERIES}\\
\noindent{\em Feature}: Computes the Poincare series for the algebras of   joint covariants for the binary forms of degrees $d_1, d_2, \ldots,d_n.$\\
\noindent{\em Calling sequence}:
{\tt
COVARIANTS\_SERIES$\left([d_1, d_2, \ldots,d_n] \right);$
}\\
\noindent{\em Parameters}:\\
{\raggedright
\begin{tabular}{lcl}
$ [d_1, d_2, \ldots,d_n]$ & - & a list of degrees of $n$ binary forms.\\
$ n $ & - & an integer, $n\geq 1.$\\

\end{tabular}

\vspace{5mm}
{\noindent
{\em Command name}: {\tt KERNEL\_SERIES}\\
\noindent{\em Feature}: Computes the Poincare series for the kernel  of   Weitzenb\"ock derivation defined by $n$ Jordan block  of sizes  $d_1+1,$ $ d_2+1,$ $ \ldots,d_n.$\\
\noindent{\em Calling sequence}:
{\tt
KERNEL\_SERIES$\left([d_1, d_2, \ldots,d_n] \right);$
}\\
\noindent{\em Parameters}:\\
{\raggedright
\begin{tabular}{lcl}
$ [d_1, d_2, \ldots,d_n]$ & - & a list of sizes of the $n$ Jordan blocks.\\
$ n $ & - & an integer, $n\geq 1.$\\

\end{tabular}

\vspace{5mm}
{\noindent
{\em Command name}: {\tt BIVARIATE\_SERIES}\\
\noindent{\em Feature}: Computes the bivariate Poincare series for  the algebra of covariants of   binary form  of degree $d.$ Also, computes the bivariate Poincare series for  the kernel of the basic Weitzenb\"ock derivation.\\
\noindent{\em Calling sequence}:
{\tt
BIVARIATE\_SERIES$([d]);$
}\\
\noindent{\em Parameters}:\\
{\raggedright
\begin{tabular}{lcl}
$ d$ & - & the degree of binary form.\\

\end{tabular}

\vspace{5mm}
{\noindent
{\em Command name}: {\tt MULTIVAR\_COVARIANTS}\\
\noindent{\em Feature}: Computes the multivariate Poincar\'e series for the algebra of   joint covariants for $n$ binary forms of degrees $d_1, d_2, \ldots,d_n.$\\
\noindent{\em Calling sequence}:
{\tt
MULTIVAR\_COVARIANTS$\left([d_1, d_2, \ldots,d_n] \right);$
}\\
\noindent{\em Parameters}:\\
{\raggedright
\begin{tabular}{lcl}
$ [d_1, d_2, \ldots,d_n]$ & - & a list of degrees of $n$ binary forms.\\
$ n $ & - & an integer, $n\geq 1.$\\

\end{tabular}

\vspace{5mm}
{\noindent
{\em Command name}: {\tt MULTIVAR\_INVARIANTS}\\
\noindent{\em Feature}: Computes the multivariate Poincar\'e series for the algebra of   joint invariants for $n$ binary forms of degrees $d_1, d_2, \ldots,d_n.$\\
\noindent{\em Calling sequence}:
{\tt
MULTIVAR\_INVARIANTS$\left([d_1, d_2, \ldots,d_n] \right);$
}\\
\noindent{\em Parameters}:\\
{\raggedright
\begin{tabular}{lcl}
$ [d_1, d_2, \ldots,d_n]$ & - & a list of degrees of $n$ binary forms.\\
$ n $ & - & an integer, $n\geq 1.$\\

\end{tabular}

\section{Examples}
\subsection{Compute $\mathcal{P}(\mathcal{I}_6,z)$} Use the command
\vspace{5mm}

{\tt > INVARIANTS\_SERIES([6]);}\\
$$
\displaystyle {\frac {{z}^{8}+{z}^{7}-{z}^{5}-{z}^{4}-{z}^{3}+z+1}{ \left( {z}^{6}+{
z}^{5}+{z}^{4}-{z}^{2}-z-1 \right)  \left( {z}^{6}+{z}^{5}-z-1
 \right)  \left( -1+{z}^{2} \right)  \left( -1+z \right) }}
$$
\subsection{Compute $\mathcal{P}(\mathcal{C}_6,z)$} Use the command
\vspace{5mm}

{\tt > COVARIANTS\_SERIES([6]);}\\
$$
\displaystyle {\frac {{z}^{10}+{z}^{8}+3\,{z}^{7}+4\,{z}^{6}+4\,{z}^{5}+4\,{z}^{4}+3
\,{z}^{3}+{z}^{2}+1}{ \left( {z}^{6}+{z}^{5}+{z}^{4}-{z}^{2}-z-1
 \right)  \left( {z}^{6}+{z}^{5}-z-1 \right)  \left( -1+{z}^{2}
 \right)  \left( -1+z \right) ^{3}}}
$$

\subsection{Compute $\mathcal{P}(\mathcal{I}_{(1,2,3)},z)$} Use the command

\vspace{5mm}

{\tt > INVARIANTS\_SERIES([1,2,3]);}\\
$$
\displaystyle {\frac {{z}^{12}+{z}^{9}+2\,{z}^{8}+3\,{z}^{7}+3\,{z}^{6}+3\,{z}^{5}+2
\,{z}^{4}+{z}^{3}+1}{ \left( -1+{z}^{4} \right) ^{2} \left( -1+{z}^{3}
 \right) ^{2} \left( -1+z \right)  \left( -1+{z}^{2} \right)  \left( {
z}^{4}+{z}^{3}+{z}^{2}+z+1 \right) }}
$$

\subsection{Compute $\mathcal{P}(\mathcal{C}_{(2,2,2)},z)$} Use the command

\vspace{5mm}

{\tt > COVARIANTS\_SERIES([2,2,2]);}\\
$$
\displaystyle {\frac {{z}^{4}+4\,{z}^{2}+1}{ \left( -1+z \right) ^{3} \left( -1+{z}^
{2} \right) ^{5}}}
$$

\subsection{Compute $\mathcal{P}( \ker \mathcal{D}_{(4)},z)$} Use the command

\vspace{5mm}

{\tt > KERNEL\_SERIES([4]);}\\
$$
\displaystyle {\frac {{z}^{2}-z+1}{ \left( -1+{z}^{2} \right)  \left( -1+{z}^{3}
 \right)  \left( -1+z \right) ^{2}}}
$$

\subsection{Compute $\mathcal{P}( \ker \mathcal{D}_{(1,1,1,2)},z)$} Use the command

\vspace{5mm}

{\tt > KERNEL\_SERIES([1,1,1,2]);}\\
$$
\displaystyle {\frac {{z}^{8}+2\,{z}^{7}+7\,{z}^{6}+11\,{z}^{5}+11\,{z}^{4}+11\,{z}^
{3}+7\,{z}^{2}+2\,z+1}{ \left( -1+{z}^{2} \right) ^{3} \left( -1+{z}^{
3} \right) ^{3} \left( -1+z \right) ^{2}}}
$$
\subsection{Compute $\mathcal{P}(\mathcal{C}_4,z,t)$} Use the command
\vspace{5mm}

{\tt > BIVARIATE\_SERIES([4]);}\\
$$
\displaystyle {\frac {{t}^{4}{z}^{2}-z{t}^{2}+1}{ \left( -1+z{t}^{2} \right) 
 \left( -1+z{t}^{4} \right)  \left( -1+{z}^{2} \right)  \left( -1+{z}^
{3} \right) }}
$$

\subsection{Compute $\mathcal{P}(\mathcal{C}_{(1,1,2)},z_1,z_2,z_3,t)$} Use the command

\vspace{5mm}

{\tt > \textcolor{red}{dd:=[1,1,2]:MULTIVAR\_COVARIANTS(dd);}}\\
$$\textcolor{blue}{
{\frac {{z_{{2}}}^{2}{z_{{1}}}^{2}{z_{{3}}}^{2}{t}^{2}+tz_{{3}}{z_{{2}
}}^{2}z_{{1}}-tz_{{3}}z_{{2}}-z_{{2}}z_{{1}}z_{{3}}+z_{{2}}z_{{1}}{t}^
{2}z_{{3}}+tz_{{3}}{z_{{1}}}^{2}z_{{2}}-z_{{1}}tz_{{3}}-1}{ \left( -1+
z_{{3}}{t}^{2} \right)  \left( -1+{z_{{3}}}^{2} \right)  \left( -1+tz_
{{2}} \right)  \left( -1+z_{{3}}{z_{{2}}}^{2} \right)  \left( -1+z_{{1
}}t \right)  \left( -1+{z_{{1}}}^{2}z_{{3}} \right)  \left( -1+z_{{2}}
z_{{1}} \right) }}
}
$$

\subsection{Compute $\mathcal{P}(\mathcal{I}_{(4,4)},z_1,z_2,t)$} Use the command

\vspace{5mm}

{\tt > \textcolor{red}{dd:=[4,4]:MULTIVAR\_INVARIANTS(dd);}}\\
$$\textcolor{blue}{
-{\frac {{z_{{1}}}^{4}{z_{{2}}}^{4}+{z_{{2}}}^{2}{z_{{1}}}^{2}+1}{
 \left( -1+{z_{{2}}}^{2} \right)  \left( -1+{z_{{2}}}^{3} \right) 
 \left( -1+{z_{{1}}}^{2}z_{{2}} \right)  \left( -1+z_{{1}}z_{{2}}
 \right)  \left( -1+z_{{1}}{z_{{2}}}^{2} \right)  \left( -1+{z_{{1}}}^
{2} \right)  \left( -1+{z_{{1}}}^{3} \right) }}
}
$$

\end{document}